\documentclass[10pt,a4paper,twoside,leqno]{article}
\usepackage{amsfonts,amssymb}
\usepackage{eufrak}
\usepackage{amscd}

\newtheorem{defn}{Definition}

\vspace{0.5cm}

 4 \font\ebf=cmbx8
\font\erm=cmr8 

\usepackage{graphicx}

\begin{document}

\thispagestyle{empty}

\noindent {\bf $\psi$-Pascal and $\hat q_{\psi}$-Pascal matrices -
an accessible factory of one source identities and resulting
applications}

\vspace{0.7cm} {\it Andrzej K. Kwa\'sniewski}

\vspace{0.7cm}

{\erm Higher School of Mathematics and Applied Informatics}

\noindent{\erm  Kamienna 17, PL-15-021 Bia\l ystok , Poland}

\vspace{0.5cm}

\noindent {\ebf Summary}

{\small Recently the author  proposed two extensions of Pascal and
$q$-Pascal matrices defined here also - in the spirit  of the Ward
"Calculus of sequences" [1] promoted in the framework  of the
$\psi$- Finite Operator Calculus [2,3] . Specifications to
$q$-calculus case and Fibonomial calculus case are made explicit
as an example of abundance of new possibilities being opened. In
broader context the $\psi$-Pascal $P_{\psi}[x]$ and $\hat
q_{\psi}$-Pascal $P_{\hat q_{\psi}}[x]$ matrices appear to be as
natural as standard Pascal matrix $P[x]$ already is known to be
[4]. Among others these are a one source factory of streams of
identities and indicated resulting applications.

\vspace{0.7cm}

\section{I. On  the usage of references }
  The papers of main reference are: [1-3].
We shall take here notation from [2,3] (see below) and the results
from [1] as well as from [2,3] - for granted. For other respective
references see:  [2,3].The acquaintance with "The matrices of
Pascal and other greats" [4] is desirable. Further relevant
references of the present author are: [5] on extended finite
operator calculus of Rota and quantum groups and other [6-7]. The
reference to q-Pascal matrix is [8] Further Pascal matrix
references  for further readings are [9-14]. One may track down
there among others relations: the Pascal Matrix versus Classical
Polynomials. The book [15] is recommended and the recent reference
[16] is useful for further applications. Very recent $\psi$-Pascal
matrix reference is [17] and  also recent  further Pascal matrices
references ( far more not complete list of them ) are to be found
in [18-21] . The book of Kassel Christian [22] - makes an
intriguing link to the advanced world of related mathematics.

Before to proceed we anyhow explain -for the reader convenience -
some of the very basic of the intuitively  useful $\psi$-notation
promoted by the author [2,3,5,6]. Here $\psi$ denotes an extension
of $\langle\frac{1}{n!}\rangle_{n\geq 0}$ sequence to quite
arbitrary one (the so called - admissible) and the specific
choices are for example: Fibonomialy-extended
$\langle\frac{1}{F_n!}\rangle_{n\geq 0}$  (here $\langle F_n
\rangle$ denotes the  Fibonacci sequence ) or Gauss $q$-extended
$\langle\frac{1}{n_q!}\rangle_{n\geq 0}$ admissible sequences of
extended umbral operator calculus or just "the usual"
$\langle\frac{1}{n!}\rangle_{n\geq 0}$ common choice. We get used
to write these $q$ - Gauss and other extensions in mnemonic
convenient upside down notation [2,3,5,6]
\begin{equation}\label{eq1}
\psi_n\equiv n_\psi , x_{\psi}\equiv \psi(x)\equiv\psi_x ,
 n_\psi!=n_\psi(n-1)_\psi!, n>0 ,
\end{equation}
\begin{equation}\label{eq2}
x_{\psi}^{\underline{k}}=x_{\psi}(x-1)_\psi(x-2)_{\psi}...(x-k+1)_{\psi}
\end{equation}
\begin{equation}\label{eq3}
x_{\psi}(x-1)_{\psi}...(x-k+1)_{\psi}= \psi(x)
\psi(x-1)...\psi(x-k-1) .
\end{equation}

\vspace{2mm}

The corresponding $\psi$-binomial symbol and $\partial _{\psi}$
difference linear operator on F[[x]] (F - any field of zero
characteristics) are below  defined accordingly where following
Roman [3,3,5,6] we shall call  $\psi=\{\psi _{n}(q) \}_{n\geq 0}$;
$\psi _{n}(q) \neq 0$; $n \geq 0$ and $\psi _{0}(q)=1$ an {\it
admissible sequence}.

\vspace{2mm}

\begin{defn}
The  $\psi$-binomial symbol is defined  as follows: \\
$$\Big({n \atop k}\Big)_{\psi}=
\frac{n_{\psi}!}{k_{\psi}!(n-k)_{\psi}!}=\Big({n \atop {n-k}}\Big
)_{\psi}$$
\end{defn}
\begin{defn}
Let  $\psi$  be admissible. Let $\partial _{\psi}$ be the linear
operator lowering degree of polynomials by one defined according
to $\partial _{\psi}x^{n}=n_{\psi}x^{n-1}$;  $n \geq 0$. Then
$\partial _{\psi}$ is called the {\it$\psi$-derivative.}
\end{defn}
You may consult [2,3,5,6]  and references therein for further
development and use of this notation  "$q$-commuting variables" -
included.

\vspace{2mm}

\section{II. Towards $\psi$-Pascal matrix factory of identities}

\vspace{2mm}

Let us define analogously to  [4,9,10] define the $\psi$-Pascal
matrix as

$$
P_{\psi}[x] = exp_{\psi}\{xK_{\psi}\}
$$
where ($Z_n$ denotes the additive cyclic group)
$$
K_{\psi} =\Big ((j+1)_{\psi}\delta_{i,j+1}\Big )_{i,j\in Z_n}
$$
therefore
$$
P_{\psi}[x] = \Big (x^{i-j}\Big({i \atop j}\Big )_{\psi}\Big
)_{i,j\in Z_n}
$$ due to:   $\partial _\psi P_{\psi}[x] =
K_{\psi}\psi P_{\psi}[x] $ where  $ \psi P_{\psi}[x]|_{x=0} =
K_{\psi} . $

Explicitly (see [8] for $q$-case) $K_{\psi}$ matrix is of the form

\vspace{2mm}

$$ K_{\psi}=\left[\begin{array}{cccccccccccccccc}
 0 & 0 & 0 & 0 & 0 & 0 & 0 & 0 & 0 & 0 & 0 & 0 & 0 & 0 & 0 & 0\\
1_{\psi}& 0 & 0 & 0 & 0 & 0 & 0 & 0 & 0 & 0 & 0 & 0 & 0 & 0 & 0 & 0\\
0 & 2_{\psi} & 0 & 0 & 0 & 0 & 0 & 0 & 0 & 0 & 0 & 0 & 0 & 0 & 0 & 0\\
0 & 0 & 3_{\psi} & 0 & 0 & 0 & 0 & 0 & 0 & 0 & 0 & 0 & 0 & 0 & 0 & 0\\
0 & 0 & 0 & 4_{\psi} & 0 & 0 & 0 & 0 & 0 & 0 & 0 & 0 & 0 & 0 & 0 & 0\\
0 & 0 & 0 & 0 & 5_{\psi} & 0 & 0 & 0 & 0 & 0 & 0 & 0 & 0 & 0 & 0 & 0\\
0 & 0 & 0 & 0 & 0 & 7_{\psi} & 0 & 0 & 0 & 0 & 0 & 0 & 0 & 0 & 0 & 0\\
0 & 0 & 0 & 0 & 0 & 0 & 8_{\psi} & 0 & 0 & 0 & 0 & 0 & 0 & 0 & 0 & 0\\
0 & 0 & 0 & 0 & 0 & 0 & 0 & 9_{\psi}& 0 & 0 & 0 & 0 & 0 & 0 & 0 & 0\\
0 & 0 & 0 & 0 & 0 & 0 & 0 & 0 & 0 & 10_{\psi} & 0 & 0 & 0 & 0 & 0 & 0\\
0 & 0 & 0 & 0 & 0 & 0 & 0 & 0 & 0 & 0 & 11_{\psi} & 0 & 0 & 0 & 0 & 0\\
 . & . & . & . & . & . & . & . & . & . & . & . & . & . & . & .\\
 . & . & . & . & . & . & . & . & . & . & . & . & . & . & . & .\\
 . & . & . & . & . & . & . & . & . & . & . & . & . & . & . & .\\
 . & . & . & . & . & . & . & . & . & . & . & . & . & . & . & .\\
 0 & 0 & 0 & 0 & 0 & 0 & 0 & 0 & 0 & 0 & 0 & 0 & 0 & 0 & (n-1)_{\psi}& 0\\
 \end{array}\right]$$

$$\textbf{ Fig.1. The $K_{\psi}$ matrix }$$

 \vspace{2mm}

Naturally    $K_{\psi}^n = 0 ; K_{\psi}^k\neq 0$ for $0 \leq k
\leq (n-1). $ Hence we have
$$
 P_{\psi}[x] =  exp_{\psi}\{xK_{\psi}\} = \sum_{k \in
Z_n} \frac{x^k K_{\psi}^k}{k_{\psi}!}
$$
the result $P_{\psi}[x]$ of $\psi$- exponentiation above  being
shown on the Fig.2.

 \vspace{2mm}

$$P_{\psi}[x]=\left[\begin{array}{ccccccccc}
 1 & 0 & 0 & 0 & 0 & 0 & 0 & 0 & 0\\
x^1 & 1 & 0 & 0 & 0 & 0 & 0 & 0 & 0\\
x^2 & 2_{\psi}x & 1 & 0 & 0 & 0 & 0 & 0 & 0\\
x^3 & 3_{\psi}x^2 & 3_{\psi}x & 1 & 0 & 0 & 0 & 0 & 0\\
x^4 & 4_{\psi}x^3 & 6_{\psi}x^2 & 4_{\psi}x & 1 & 0 & 0 & 0 & 0\\
. & . & . & . & . & . & . & . & 0\\
. & . & . & . & . & . & . & . & 0\\
. & . & . & . & . & . & . & . & 0\\
x^{n-1} & 0 & 0 & 0 & 0 & 0 & 0 & (n-1)_{\psi}x & 1\\
 \end{array}\right]$$

$$\textbf{ Fig. 2. The $P_{\psi}[x]$ matrix }$$

\vspace{2mm}

Immediately we see that  the $\psi$-Pascal matrix $P_{\psi}[x]
=exp_{\psi}\{xK_{\psi}\}$ is also the source of many important
identities. Here below there are the examples correspondent to
those from [4] which are accordingly infered from the
$\psi$-additivity property (non-group property in general) :

$$
P_{\psi}[x]P_{\psi}[y] = P_{\psi}[x+_{\psi} y].
$$

 \vspace{2mm}

\textbf{Warning}: for not normal sequences : see: [1,2,3,5,6,8] -
the one parameter family $\{P_{\psi}[x]\}_{x \in F}$ is \textit{
not a group }! since for not normal sequences $(1-_{\psi}
1)^{2k}\neq 0$ thought $[x +_{\psi} (-x) ]^{2k+1} = 0.$

 \vspace{2mm}

In general we are dealing with abelian semigroup with identity
which becomes the group only for normal sequences. And so coming
back to identities we have for example :

 \vspace{2mm}

\begin{equation}\label{eq4}
\sum_{j \leq k \leq i}\Big({i \atop k}\Big )_{\psi}\Big({k \atop
j}\Big )_{\psi}= (1 +_{\psi}1)^{i-j}\Big({i \atop j}\Big )_{\psi}
,  i\geq j  \Longleftrightarrow P_{\psi}[1] P_{\psi}[1] =
P_{\psi}[1 +_{\psi} 1] .
\end{equation}

\begin{equation}\label{eq5}
\sum_{j \leq k \leq i}(-1)^k\Big({i \atop k}\Big )_{\psi}\Big({k
\atop j}\Big )_{\psi}= (1 -_{\psi}1)^{i-j}\Big({i \atop j}\Big
)_{\psi} ,  i\geq j  \Longleftrightarrow P_{\psi}[1] P_{\psi}[-1]
= P_{\psi}[1 -_{\psi} 1] .
\end{equation}

\vspace{2mm}

The above identities after the choice $\psi =
\langle\frac{1}{n!}\rangle_{n\geq 0}$ coincide with the
corresponding ones from [4]. There are much more examples of this
nature.

\vspace{2mm}

We shall now try also to find out a kind of $\psi$-extended
version of the $q$-identity (6)

\begin{equation}\label{eq6}
\sum_{0 \leq k \leq i} \Big({i \atop k}\Big)^2_q = \Big({2i \atop
i}\Big )_q \Longleftrightarrow P_{\psi}[1] P^T_{\psi}[1] \equiv
F_q[1] .
\end{equation}

\vspace{2mm}

where  we have defined the $q$-Fermat matrix as follows

\begin{equation}\label{eq7}
F_q[1] =  \Big(\Big({i+j \atop i}\Big )_q \Big )_{i,j\in Z_n}.
\end{equation}

\vspace{2mm}

For q =1 case- name Fermat -  see [15] for this Fermat called
Pascal symmetric Matrix for q =1 see:  [4,9]. For $q$-binomial -
see below in \textbf{Important}.

\vspace{2mm}

In order to find out a kind of  $\psi$-extended version of the
Pascal-Fermat $q$-identity identity (6)  we shall proceed as in
[16].  There the Cauchy $\hat {q}_{\psi} $- identity and $\hat
{q}_{\psi}$-Fermat matrix were introduced  due to the use of the
$\hat {q}_{\psi} $-muting variables from Extended Finite Operator
Calculus [3,5]. The \textit{linear}  $\hat {q}_{\psi}$-mutator
operator  was  defined in [3,5,16]  as follows for $F$ - field of
characteristic zero and $F[x]$ - the linear space of polynomials.

\vspace{2mm}
$$
 \hat {q}_{\psi} :F[x] \to F[x] ;
\quad \hat {q}_{\psi} x^n = \frac{{\left( {n + 1} \right)_{\psi}
- 1}}{{n_{\psi} } }x^n ;\quad n \ge 0.
$$
\vspace{2mm}

\textbf{Important}. With the Gaussian choice of admissible
sequence [3,5]
$$
\psi= \{\psi_n(q)\}_{n\geq 0},   \psi_n(q)=[n_q!]^{-1} ,
n_q=\frac{1-q^n}{1-q} ,  n_q!=n_q(n-1)_q! , 1_q!=0_q!=1 , \hat
{q}_{\psi}x^n = q^nx^n
$$
and  the $\hat {q}_{\psi}$-Pascal and $\hat {q}_{\psi}$-Fermat
matrices from [16] (see next section ) coincide with $q$-Pascal
and $q$-Fermat matrices correspondingly \textbf{which is not the
case} for the general case - for example Fibonomial $F$-Pascal
matrix is different from $\hat {q}_F$-Pascal matrix - see next
section.

\vspace{2mm}

In [16] in analogy to the standard case [9,10,4] the matrices with
operator valued matrix elements

$$
x^{i-j}\Big({i \atop j}\Big )_{\hat q_{\psi}} , \Big({{i+j} \atop
j}\Big )_{\hat q_{\psi}}
$$
were named the  $\hat q_{\psi}$- Pascal P[x]  and  $\hat q_{\psi}$
-Fermat F[1]  matrices - correspondingly  i.e.

\vspace{2mm}
$$
P_{\hat q_{\psi}}[x] = \Big (x^{i-j}\Big({i \atop j}\Big )_{\hat
q_{\psi}}\Big )_{i,j\in Z_n}
$$

\vspace{2mm}

The $\hat {q}_{\psi}$-P[1] Pascal and $\hat {q}_{\psi}$-F[1]
Fermat matrices from [16] are related via the following identity
for operator valued matrix elements

\vspace{2mm}

\begin{equation}\label{}
\sum_{k\geq 0} {\hat q_{\psi}}^{(r-k)(j-k)}\Big({i \atop k}\Big
)_{\hat q_{\psi}}\Big({j\atop k}\Big )_{\hat q_{\psi}} =
\Big({{i+j} \atop j}\Big )_{\hat q_{\psi}} .
\end{equation}

\vspace{2mm}

The relation (8)  is the one being looked for to extend the
Pascal-Fermat $q$-identity  (6). Here -  following [16]- we use
the new  $\hat q_{\psi}$-Gaussian symbol with operator valued
matrix elements.

\vspace{2mm}

\begin{defn}
We define $\hat{q}_{\psi}$-binomial symbol i.e. $\hat
{q}_{\psi}$-Gaussian coefficients as follows: \\
$\Big({n \atop k}\Big)_{\hat {q}_{\psi}}=
\frac{n_{\hat{q}_{\psi}}!}{k_{\hat
q_{\psi}}!(n-k)_{\hat{q}_{\psi}}!}=\Big({n \atop {n-k}}\Big
)_{\hat {q}_{\psi}}$ { where } $n_{\hat {q}_{\psi}}!=n_{\hat
{q}_{\psi}}(n-1)_{\hat {q}_{\psi}}! , 1_{\hat {q}_{\psi}}!=0_{\hat
{q}_{\psi}}!=1$ and $n_{\hat {q}_{\psi}}=\frac{1-{\hat
{q}_{\psi}}^n}{1-{\hat {q}_{\psi}}}$ for $n>0$.
\end{defn}

\section{III. Specifications : $q$-umbral and umbral Fibonomial cases}

\vspace{2mm}

\textbf{III-q} $q$-umbral calculus case [1,2,3,5-8]

\vspace{2mm}

Let us make the $q$-Gaussian choice [2,3,5,6,8]  of the admissible
sequence $\psi= \{\psi_n(q)\}_{n\geq 0}$. Then the $\psi$-Pascal
matrix  becomes the $q$-Pascal matrix from [8] and we arrive
mnemonic at the corresponding to $q=1$ case numerous
$q$-identities and other "$q$-applications". Specifically in the
$q$-case  we have (see: Proposition 4.2.3  in [22])

\vspace{2mm}

\begin{equation}\label{}
\sum_{k\geq 0}  q^{(r-k)(j-k)}\Big({r \atop k}\Big )_q\Big({s
\atop{j-k}}\Big )_q = \Big({{r+s} \atop j}\Big )_q
\end{equation}

\vspace{2mm}

hence from this  Cauchy $q$-identity we obtain the following easy
to find out formula for the symmetric Pascal (or Fermat) matrix
elements:

\vspace{2mm}

\begin{equation}\label{}
\sum_{k\geq 0} q^{(r-k)(j-k)}\Big({i \atop k}\Big )_q \Big({j\atop
k}\Big )_q = \Big({{i+j} \atop j}\Big )_q .
\end{equation}

\vspace{2mm}

Naturally we are dealing now with not normal sequences i.e. not
with a one parameter $q$-Pascal group [8]  since for  $(1-_q
1)^{2k} \neq 0$ though $[x +_q (-x)]^{2k+1} = 0$ ;  see: [1] and
then [2,3,5,6,8]. If $q$-Pascal matrix $P_q[1]=
exp_q\{xK_q\}|_{x=1}$ is considered also for $q \in {GF(q)}$ field
then $q = p_m$ where $p$ is prime  and  $\Big({n \atop k}\Big )_q$
becomes the number of $k$-dimensional subspaces in $n-th$
dimensional space over Galois field  $GF(q)$. Also  $q$  real and
$-1<q<+1$ cases are exploited in vast literature on the so-called
$q$-umbral calculus (for Cigler , Roman and Others see: [3,23] and
references therein- links to thousands in [23]). It is not
difficult to notice that the $\hat q_\psi $- Pascal and $\hat
q_\psi $-Fermat matrices under the $q$-Gassian choice of the
admissible sequence $\psi$ -  coincide with $q$-Pascal and
$q$-Fermat matrices correspondingly which is meaningful
\textbf{magnificent exception} and which is not the case in
general .

\vspace{2mm}

\textbf{III-F} \textbf{FFOC}-umbral calculus case [6-7]

\vspace{2mm}

 In straightforward analogy to the $q$-case above
consider now the Fibonomial coefficients  (see:  FFOC = Fibonomial
Finite Operator Calculus  Example 2.1 in [6]) where  $F_n$ denote
the Fibonacci numbers and    $\psi_n(q)=[F_n!]^{-1}$ .
$$
\left( \begin{array}{c} n\\k\end{array}
\right)_{F}=\frac{F_{n}!}{F_{k}!F_{n-k}!}\equiv
\frac{n_{F}^{\underline{k}}}{k_{F}!},\quad n_{F}\equiv F_{n}\neq
0, $$

\noindent where we make an analogy driven [6,5,3,2]
identifications $(n>0)$:
$$
n_{F}!\equiv n_{F}(n-1)_{F}(n-2)_{F}(n-3)_{F}\ldots 2_{F}1_{F};$$
$$0_{F}!=1;\quad n_{F}^{\underline{k}}=n_{F}(n-1)_{F}\ldots (n-k+1)_{F}. $$

\textbf{Information}
 In [7]  a partial ordered set was defined in such a way that
 the Fibonomial coefficients count the number of specific finite
"birth-self-similar"  sub-posets of this  infinite non-tree poset
naturally related to the Fibonacci tree of rabbits growth process.

\vspace{2mm}

 The $\psi$-Pascal matrix  becomes then the $F$-Pascal
matrix and we arrive at the corresponding $F$-identities (mnemonic
replacement of  $\psi$  by  $F$ ) and other "$F$-applications" -
hoped to be explored soon.

\vspace{2mm}

Naturally we are now dealing with \textbf{not normal} sequences :
see: [1,2,3,5,6,8] -i.e.we have no  $F$-Pascal \textbf{group}
since for $(1-_F 1)^{2k} \neq 0$  though   $[x +_F (-x)]^{2k+1} =
0$. For example: $(x +_F y)^2 = x^2 +  F_2xy + y^2,   (x +_F y)^4
= x^4 + F_4 x^3y + F_4F_3x^2y^2 + F_4 xy^3 + y^4$.

\vspace{2mm}

Here in the Fibonomial choice case the semi-group generating
matrix  matrix $K_F$ is of the form

\vspace{2mm}

\vspace{2mm}

$$K_F = \left[\begin{array}{cccccccccccccccc}
 0 & 0 & 0 & 0 & 0 & 0 & 0 & 0 & 0 & 0 & 0 & 0 & 0 & 0 & 0 & 0\\
F_1& 0 & 0 & 0 & 0 & 0 & 0 & 0 & 0 & 0 & 0 & 0 & 0 & 0 & 0 & 0\\
0 & F_2 & 0 & 0 & 0 & 0 & 0 & 0 & 0 & 0 & 0 & 0 & 0 & 0 & 0 & 0\\
0 & 0 & F_3 & 0 & 0 & 0 & 0 & 0 & 0 & 0 & 0 & 0 & 0 & 0 & 0 & 0\\
0 & 0 & 0 & F_4 & 0 & 0 & 0 & 0 & 0 & 0 & 0 & 0 & 0 & 0 & 0 & 0\\
0 & 0 & 0 & 0 & F_5 & 0 & 0 & 0 & 0 & 0 & 0 & 0 & 0 & 0 & 0 & 0\\
0 & 0 & 0 & 0 & 0 & F_7 & 0 & 0 & 0 & 0 & 0 & 0 & 0 & 0 & 0 & 0\\
0 & 0 & 0 & 0 & 0 & 0 & F_8 & 0 & 0 & 0 & 0 & 0 & 0 & 0 & 0 & 0\\
0 & 0 & 0 & 0 & 0 & 0 & 0 & F_9 & 0 & 0 & 0 & 0 & 0 & 0 & 0 & 0\\
0 & 0 & 0 & 0 & 0 & 0 & 0 & 0 & F_{10} & 0 & 0 & 0 & 0 & 0 & 0 & 0\\
0 & 0 & 0 & 0 & 0 & 0 & 0 & 0 & 0 & F_{11} & 0 & 0 & 0 & 0 & 0 & 0\\
 . & . & . & . & . & . & . & . & . & . & . & . & . & . & . & .\\
 . & . & . & . & . & . & . & . & . & . & . & . & . & . & . & .\\
 . & . & . & . & . & . & . & . & . & . & . & . & . & . & . & .\\
 . & . & . & . & . & . & . & . & . & . & . & . & . & . & . & .\\
 0 & 0 & 0 & 0 & 0 & 0 & 0 & 0 & 0 & 0 & 0 & 0 & 0 & 0 & F_{(n-1)}& 0\\
 \end{array}\right]$$

$$\textbf{ Fig.3. The $K_F$ matrix }$$

and the corresponding beautiful $F$-Pascal matrix $P_F[x] = exp_F
\{xK_F\}$ being $F$-exponentiation of  $K_F$ reads:

 \vspace{2mm}

$$ P_F[x]= \left[\begin{array}{ccccccccc}
 1 & 0 & 0 & 0 & 0 & 0 & 0 & 0 & 0\\
x^1 & 1 & 0 & 0 & 0 & 0 & 0 & 0 & 0\\
x^2 & F_2x & 1 & 0 & 0 & 0 & 0 & 0 & 0\\
x^3 & F_3x^2 & F_3x & 1 & 0 & 0 & 0 & 0 & 0\\
x^4 & F_4x^3 & F_6x^2 & F_4x & 1 & 0 & 0 & 0 & 0\\
. & . & . & . & . & . & . & . & 0\\
. & . & . & . & . & . & . & . & 0\\
. & . & . & . & . & . & . & . & 0\\
x^{n-1} & F_{n-1}x^{n-1} & \Big({{n-1} \atop 2}\Big )_Fx^{n-2} & . & . & . & . & F_{(n-1)}x & 1\\
 \end{array}\right]$$

$$\textbf{ Fig. 4. The $P_F[x]$ matrix }$$

\vspace{2mm}

The $\hat q_F$-Pascal $P_{K_{\hat q_F}}[ x]$  and $K_{\hat
q_F}$-Fermat matrix do not coincide with $F$-Pascal and $F$-Fermat
matrices correspondingly  as indicated earlier though in our
friendly-mnemonic notation they look so much alike . Namely , the
corresponding $K_{\hat q_{\psi}}$ matrix with the Fibonomial
choice  $\psi_n(q)=[F_n!]^{-1}$ is now of the form

\vspace{2mm}

$$K_{\hat q_F}= \left[\begin{array}{cccccccccccccccc}
 0 & 0 & 0 & 0 & 0 & 0 & 0 & 0 & 0 & 0 & 0 & 0 & 0 & 0 & 0 & 0\\
1_{\hat q_F}& 0 & 0 & 0 & 0 & 0 & 0 & 0 & 0 & 0 & 0 & 0 & 0 & 0 & 0 & 0\\
0 & 2_{\hat q_F} & 0 & 0 & 0 & 0 & 0 & 0 & 0 & 0 & 0 & 0 & 0 & 0 & 0 & 0\\
0 & 0 & 3_{\hat q_F} & 0 & 0 & 0 & 0 & 0 & 0 & 0 & 0 & 0 & 0 & 0 & 0 & 0\\
0 & 0 & 0 & 4_{\hat q_F} & 0 & 0 & 0 & 0 & 0 & 0 & 0 & 0 & 0 & 0 & 0 & 0\\
0 & 0 & 0 & 0 & 5_{\hat q_F} & 0 & 0 & 0 & 0 & 0 & 0 & 0 & 0 & 0 & 0 & 0\\
0 & 0 & 0 & 0 & 0 & 7_{\hat q_F} & 0 & 0 & 0 & 0 & 0 & 0 & 0 & 0 & 0 & 0\\
0 & 0 & 0 & 0 & 0 & 0 & 8_{\hat q_F} & 0 & 0 & 0 & 0 & 0 & 0 & 0 & 0 & 0\\
0 & 0 & 0 & 0 & 0 & 0 & 0 & 9_{\hat q_F}& 0 & 0 & 0 & 0 & 0 & 0 & 0 & 0\\
0 & 0 & 0 & 0 & 0 & 0 & 0 & 0 & 0 & 10_{\hat q_F} & 0 & 0 & 0 & 0 & 0 & 0\\
0 & 0 & 0 & 0 & 0 & 0 & 0 & 0 & 0 & 0 & 11_{\hat q_F} & 0 & 0 & 0 & 0 & 0\\
 . & . & . & . & . & . & . & . & . & . & . & . & . & . & . & .\\
 . & . & . & . & . & . & . & . & . & . & . & . & . & . & . & .\\
 . & . & . & . & . & . & . & . & . & . & . & . & . & . & . & .\\
 . & . & . & . & . & . & . & . & . & . & . & . & . & . & . & .\\
 0 & 0 & 0 & 0 & 0 & 0 & 0 & 0 & 0 & 0 & 0 & 0 & 0 & 0 & (n-1)_{\hat q_F}& 0\\
 \end{array}\right]$$

$$\textbf{ Fig.5. The $K_{\hat q_F}$ matrix }$$

 \vspace{2mm}

Similarly to the earlier case considered $K_{\hat q_F}^n = 0 ;
K_{\hat q_F}^k\neq 0$ for $0 \leq k \leq (n-1)$ and again we also
 have
$$
 P_{\hat q_F}[x] =  exp_{\psi}\{xK_{\hat q_F}\} = \sum_{k \in
Z_n} \frac{x^kK_{\hat q_F}^k}{k_F!}.
$$
The result $P_{\hat q_F}[x]$ of the $F$-exponentiation above has
been shown on the Fig.6.

 \vspace{2mm}

$$ P_{\hat q_F}[x]= \left[\begin{array}{ccccccccc}
 1 & 0 & 0 & 0 & 0 & 0 & 0 & 0 & 0\\
x_{\psi}& 1 & 0 & 0 & 0 & 0 & 0 & 0 & 0\\
x^2 & 2_{\psi}x & 1 & 0 & 0 & 0 & 0 & 0 & 0\\
x^3 & 3_{\psi}x^2 & 3_{\psi}x & 1 & 0 & 0 & 0 & 0 & 0\\
x^4 & 4_{\psi}x^3 & 6_{\psi}x^2 & 4_{\psi}x & 1 & 0 & 0 & 0 & 0\\
. & . & . & . & . & . & . & . & 0\\
. & . & . & . & . & . & . & . & 0\\
. & . & . & . & . & . & . & . & 0\\
x^{n-1} & 0 & 0 & 0 & 0 & 0 & 0 & (n-1)_{\psi}x & 1\\
 \end{array}\right]$$

$$\textbf{ Fig.6. The $P_{\hat q_F}[x]$ matrix }$$

\vspace{2mm}

\textbf{Important-Conclusive.}  Apart then from $\psi$-Pascal one
source matrix factory of identities we indicate in explicit also
the origins of the ${\hat q_F}$ - Pascal and ${\hat q_F}$-Fermat
matrices factory of mnemonic attainable identities (compare via
[16] with [9-14,18-21,4] ). From operator identities involving the
$\hat q_{\psi}$-Pascal $P_{K_{\hat q_{\psi}}}[ x]$ and $K_{\hat
q_{\psi}}$-Fermat matrix we obtain identities in terms of objects
on which the $\hat q_{\psi}$ (or $\hat q_{\psi,Q}$ from
[3,5,6,16]) act and these are polynomials from F[x] or in more
general setting [6,5,3] from formal series algebra F[[x]] where F
denotes any field of zero characteristics. In order to get such
countless realizations of operator identities in terms of formal
series it is enough to act by both sides of a given operator
identity on the same element from $F[[x]]$.

\section{IV.  Remark on perspectives}

\vspace{2mm}

The perspective of numerous applications are opened. Apart from
being the natural one source factory of identities $\psi$-Pascal
$P_{\psi}[x]$ and   $\hat q_{\psi}$-Pascal $P_{K_{\hat q_{\psi}}}[
x]$ and $K_{\hat q_{\psi}}$-Fermat matrices as well appear  to be
the similar way natural objects and tools as the Pascal matrix
P[x] is in the already  mentioned and other applications - (see
[4,18]- for example). Just to indicate few more of them:  the
considerations and results of [4] concerned with Bernoulli
polynomials might be extended  to the case of $\psi$-basic
Bernoulli-Ward polynomials introduced in [1] and investigated
recently in [17] in the framework of the $\psi$- Finite Operator
Calculus [2,3,5-7] due to the use of  the  $\psi$- integration
proposed in [2,6] .  The same applies equally well to the case of
$\psi$-basic Hermite-Ward polynomials and other examples of
$\psi$-basic generalized Appell polynomials [3,2,5-6] which -
being  of course $\psi$- Sheffer are characterized equivalently by
the familiar $\psi$-Sheffer identity [3,2]

\begin{equation}\label{}
A_n(x +_{\psi}y) = \sum_{k \geq 0}\Big({n \atop k}\Big)_{\psi}
A_k(y)x^{n-k}.
\end{equation}

\vspace{2mm}

For further possibilities - see references [8-14,18-21] and many
other ones not known for the moment to the present author.

\end{document}